\documentclass[12pt]{article}

\usepackage{theorem,amssymb}

\advance \topmargin by -\headheight
\advance \topmargin by -\headsep
\evensidemargin \oddsidemargin
\author{J.-P. Allouche\footnote{The author was partially
supported by the ANR project ``FAN'' (Fractal and Numeration).}\\
CNRS, Institut de Math\'ematiques de Jussieu-PRG \\
\'Equipe Combinatoire et Optimisation \\
Universit\'e Pierre et Marie Curie, Case 247 \\
4 Place Jussieu \\
F-75252 Paris Cedex 05 France \\
{\tt jean-paul.allouche@imj-prg.fr}
}

\title{A note on products involving $\zeta(3)$ and Catalan's constant}

\date{ }

\def \proof{\bigbreak\noindent{\it Proof.\ \ }}

\def \endpf{{\ \ $\Box$ \medbreak}}

\newtheorem{proposition}{Proposition}
\theorembodyfont{\rm}

\newtheorem{remark}{Remark}

\begin{document}

\maketitle

\begin{abstract}
In a recent paper Kachi and Tzermias give elementary proofs of four product formulas involving
$\zeta(3)$, $\pi$, and Catalan's constant. They indicate that they were not able to deduce 
these products directly from the values of a function introduced in 1993 by Borwein and 
Dykshoorn. We provide here such a proof for two of these formulas. We also give a direct 
proof for the other two formulas, by using a generalization of the Borwein-Dykshoorn 
function due to Adamchik. Finally we give an expression of the Borwein-Dykshoorn function
in terms of the ``parameterized-Euler-constant function'' introduced by Xia in 2013, which 
happens to be a particular case of the ``generalized Euler constant function'' introduced by
K. and T. Hessami Pilehrood in 2010.

\medskip


\medskip

\noindent
{\bf MSC Classes}: 11Y60, 33B99, 11M06, 11M99

\end{abstract}

\section{Introduction}
In a recent paper Kachi and Tzermias prove in an elementary way four nice formulas involving 
$\zeta(3)$, $\pi$, and Catalan's constant (see \cite[Propositions~1 and 2]{KT}), namely

\begin{equation}\label{1}
\lim_{n \to \infty} \prod_{k=1}^{2n+1} e^{-1/4} \left(1-\frac{1}{k+1}\right)^{\frac{k(k+1)}{2}(-1)^k}
= \exp\left(\frac{7\zeta(3)}{4\pi^2} + \frac{1}{4}\right) 
\end{equation}

\begin{equation}\label{2}
\lim_{n \to \infty} \prod_{k=1}^{2n} e^{1/4} \left(1-\frac{1}{k+1}\right)^{\frac{k(k+1)}{2}(-1)^k}
= \exp\left(\frac{7\zeta(3)}{4\pi^2} - \frac{1}{4}\right)
\end{equation}

\begin{equation}\label{3}
\lim_{n \to \infty} \prod_{k=1}^{2n} \left(1-\frac{2}{2k+1}\right)^{k(-1)^k}
= \exp\left(\frac{2G}{\pi} - \frac{1}{2}\right)
\end{equation}

\begin{equation}\label{4}
\lim_{n \to \infty} \prod_{k=1}^{2n+1} \left(1-\frac{2}{2k+1}\right)^{k(-1)^k}
= \exp\left(\frac{2G}{\pi} + \frac{1}{2}\right)
\end{equation}
where $G = \sum_{n \geq 0} \frac{(-1)^n}{(2n+1)^2}$ is the Catalan constant.

\medskip

At the end of their article \cite{KT} the authors recall the Borwein-Dykshoorn function
$$
D(x) = \lim_{n \to \infty} \prod_{k=1}^{2n+1} \left(1+\frac{x}{k}\right)^{k(-1)^{k+1}}
$$
which was introduced in \cite{BD} as a generalization of a result of Melzak \cite{Melzak}
proving that
$$
\lim_{n \to \infty} \prod_{k=1}^{2n+1} \left(1+\frac{2}{k}\right)^{k(-1)^{k+1}} = \frac{\pi e}{2}\cdot
$$
Kachi and Tzermas indicate that they were not able to deduce any of the relations~(\ref{1}), (\ref{2}),
(\ref{3}), and (\ref{4}) from \cite{BD} (where the authors give the values of $D(a/b)$ for $a$ integer 
and $b = 1, 2, 3$), though, e.g., the constant $e^{G\pi}$ occurs both in \cite{BD} and in \cite{KT}).

\medskip

In this paper we give a direct proof of relations~(\ref{3}) and (\ref{4}) using the function $D(x)$. 
(Actually we only need the values $D(1)$ and $D(\frac{1}{2})$.) We then use a function similar to $D(x)$
introduced by Adamchik in \cite[p.~284]{Adam}, namely
$$
E(x) = \lim_{n \to \infty} \prod_{k=1}^{2n} \left(1-\frac{4x^2}{k^2}\right)^{-k^2(-1)^k}
$$
to prove directly relations~(\ref{1}) and (\ref{2}). (Actually we only use the value $E(\frac{1}{2})$.)

\section{Formulas~(\ref{3}) and (\ref{4})}

\begin{proposition}
Formulas~(\ref{3}) and (\ref{4}) can be deduced directly from the values of the Borwein-Dykshoorn
function $D(1)$ and $D(1/2)$, and from classical results for the function $\Gamma$.
\end{proposition}

\proof

We first note that
$$
\prod_{k=1}^{2n+1} \left(1 - \frac{2}{2k+1}\right)^{k(-1)^k} =
\left(1 - \frac{2}{4n+3}\right)^{-(2n+1)}
\prod_{k=1}^{2n} \left(1 - \frac{2}{2k+1}\right)^{k(-1)^k}
$$
Since $\lim_{n \to \infty} \left(1 - \frac{2}{4n+3}\right)^{-(2n+1)} = e$ (take the
logarithm) it is clear that (\ref{4}) is readily deduced from (\ref{3}). It thus suffices
to prove (\ref{3}).

\medskip

Let $D_n(x) := \displaystyle\prod_{k=1}^{2n+1} \left(1 + \frac{x}{k}\right)^{k(-1)^{k+1}}$ and
$A_n := \displaystyle\prod_{k=1}^{2n} \left(1 - \frac{2}{2k+1}\right)^{k(-1)^k}$. Then
$$
\begin{array}{lll}
A_n &=&
\displaystyle\prod_{k=1}^{2n} \left(\frac{2k-1}{2k+1}\right)^{k(-1)^k} 
= \displaystyle\prod_{k=1}^{2n} \left(\frac{2k-1}{2k}\right)^{k(-1)^k} 
\displaystyle\prod_{k=1}^{2n} \left(\frac{2k+1}{2k}\right)^{k(-1)^{k+1}} \\
&=& \displaystyle\prod_{k=1}^{2n} \left(\frac{2k-1}{2k}\right)^{k(-1)^k}
D_{n-1}\left(\frac{1}{2}\right){\left(\displaystyle\frac{4n}{4n+1}\right)^{2n}} \\
\end{array}
$$
But
$$
\begin{array}{lll}
\displaystyle\prod_{k=1}^{2n} \left(\frac{2k-1}{2k}\right)^{k(-1)^k}
&=& \displaystyle\prod_{\ell=0}^{2n-1} \left(\frac{2\ell+1}{2\ell+2}\right)^{(\ell+1)(-1)^{\ell+1}} \\
&=& \displaystyle\prod_{\ell=0}^{2n-1} \left(\frac{2\ell+1}{2\ell+2}\right)^{\ell(-1)^{\ell+1}} \times
\displaystyle\prod_{\ell=0}^{2n-1} \left(\frac{2\ell+1}{2\ell+2}\right)^{(-1)^{\ell+1}} \\
&=& \displaystyle\prod_{\ell=1}^{2n-1} \left(\frac{2\ell+1}{2\ell+2}\right)^{\ell(-1)^{\ell+1}} \times
\prod_{\ell=0}^{n-1} \left(\frac{4\ell+2}{4\ell+1} \cdot \frac{4\ell+3}{4\ell+4}\right) \\
&=& \displaystyle\prod_{\ell=1}^{2n-1} \left(\frac{2\ell+1}{2\ell}\right)^{\ell(-1)^{\ell+1}} \times
\displaystyle\prod_{\ell=1}^{2n-1} \left(\frac{2\ell}{2\ell+2}\right)^{\ell(-1)^{\ell+1}} \times
\prod_{\ell=0}^{n-1} \left(\frac{4\ell+2}{4\ell+1} \cdot \frac{4\ell+3}{4\ell+4}\right) \\
&=&
\displaystyle\frac{\displaystyle\prod_{\ell=1}^{2n-1}\left(1+\frac{1}{2\ell}\right)^{\ell(-1)^{\ell+1}}}
{\displaystyle\prod_{\ell=1}^{2n-1}\left(1+\frac{1}{\ell}\right)^{\ell(-1)^{\ell+1}}} \times
\prod_{\ell=0}^{n-1} \left(\frac{4\ell+2}{4\ell+1} \cdot \frac{4\ell+3}{4\ell+4}\right) \\
&=&
\displaystyle\frac{D_{n-1}\left(\frac{1}{2}\right)}{D_{n-1}(1)} \times
\prod_{\ell=0}^{n-1} \left(\frac{4\ell+2}{4\ell+1} \cdot \frac{4\ell+3}{4\ell+4}\right) \\
\end{array}
$$
Hence
$$
A_n = 
\displaystyle
\frac{D_{n-1}\left(\frac{1}{2}\right)^2}{D_{n-1}(1)}
\left(\displaystyle\frac{4n}{4n+1}\right)^{2n} \times
\prod_{\ell=0}^{n-1} \left(\frac{4\ell+2}{4\ell+1} \cdot \frac{4\ell+3}{4\ell+4}\right)
$$
We note that
$\displaystyle\lim_{n \to \infty} \left(\frac{4n}{4n+1}\right)^{2n} = e^{-1/2}$
(take the logarithm), and that
$$
\lim_{n \to \infty}
\prod_{\ell=0}^{n-1} \left(\frac{4\ell+2}{4\ell+1} \cdot \frac{4\ell+3}{4\ell+4}\right) \\
= \lim_{n \to \infty}
\prod_{\ell=0}^{n-1} \left(\frac{\ell+1/2}{\ell+1/4} \cdot \frac{\ell+3/4}{\ell+1}\right)
= \frac{\Gamma\left(\frac{1}{4}\right)\Gamma(1)}{\Gamma\left(\frac{1}{2}\right)\Gamma\left(\frac{3}{4}\right)}
= \frac{\Gamma\left(\frac{1}{4}\right)}{\Gamma\left(\frac{1}{2}\right)\Gamma\left(\frac{3}{4}\right)}
$$
(see, e.g., \cite[Section~12-13]{WW}). Furthermore, from \cite{BD}
$$
\lim_{n \to \infty} D_{n-1}\left(\frac{1}{2}\right) = D\left(\frac{1}{2}\right) = 
\frac{2^{1/6} \sqrt{\pi} A^3 e^{G/\pi}}{\Gamma\left(\frac{1}{4}\right)}
$$
and
$$
\lim_{n \to \infty} D_{n-1}(1) = D(1) = \frac{A^6}{2^{1/6} \sqrt{\pi}}
$$
where $A$ is the Glaisher-Kinkelin constant ($A = \exp\left(\frac{1}{12} - \zeta'(-1)\right)$
where $\zeta$ is the Riemann zeta function). Putting these relations together yields
$$
\lim_{n \to \infty} A_n = \frac{\sqrt{2} \pi^{3/2}}
{\Gamma\left(\frac{1}{4}\right)\Gamma\left(\frac{3}{4}\right)\Gamma\left(\frac{1}{2}\right)}
e^{\frac{2G}{\pi} - \frac{1}{2}}.
$$
Then, Euler's reflection formula $\Gamma(z)\Gamma(1-z)= \pi/\sin(\pi z)$
(see, e.g., \cite[Section~12-14]{WW}) yields the classical relations
$\Gamma(1/2) = \sqrt{\pi}$ and $\Gamma(1/4)\Gamma(3/4) = \pi \sqrt{2}$. Thus
$$
\lim_{n \to \infty} A_n = e^{\frac{2G}{\pi} - \frac{1}{2}}
$$
which is Formula~(\ref{3}). \endpf

\section{Formulas~(\ref{1}) and (\ref{2})}

\begin{proposition}
Formulas~(\ref{1}) and (\ref{2}) can be deduced directly from the value of the Adamchik
function $E(1/2)$, and from classical results for the function $\Gamma$.
\end{proposition}

\proof  We first note that
$$
\prod_{k=1}^{2n+1} e^{-1/4} \left(1-\frac{1}{k+1}\right)^{\frac{k(k+1)}{2}(-1)^k}
= \alpha_n \left(\prod_{k=1}^{2n} e^{1/4} \left(1-\frac{1}{k+1}\right)^{\frac{k(k+1)}{2}(-1)^k}\right)
$$
with $\alpha_n = e^{-n-\frac{1}{4}} \left(\frac{2n+1}{2n+2}\right)^{-(n+1)(2n+1)}$. (It might be worth 
underlining that the first product involves $e^{-1/4}$ while the second product involves $e^{1/4}$.)
Since $\alpha_n$ tends to $e^{1/2}$ (take the logarithm), it is clear that Formula~(\ref{2}) 
implies Formula~(\ref{1}). It thus suffices to prove Formula~(\ref{2}). 

\medskip

Let $E_n = \displaystyle\prod_{k=1}^{2n} e^{1/4} \left(1 - \frac{1}{k+1}\right)^{\frac{k(k+1)}{2}(-1)^k}$.
We write $E_n^2$ in two different ways. On one hand 
$$
\begin{array}{lll}
E_n^2 &=& \displaystyle\prod_{k=1}^{2n} e^{1/2} \left(1 - \frac{1}{k+1}\right)^{(k^2+k)(-1)^k} 
       = e^n \displaystyle\prod_{k=1}^{2n}\left(1 - \frac{1}{k+1}\right)^{k^2(-1)^k}
             \times \displaystyle\prod_{k=1}^{2n}\left(1 - \frac{1}{k+1}\right)^{k(-1)^k} \\
      &=& e^n \displaystyle\prod_{k=1}^{2n}\left(1 - \frac{1}{k+1}\right)^{k^2(-1)^k}
              \times \displaystyle\prod_{\ell=2}^{2n+1}\left(1 - \frac{1}{\ell}\right)^{(\ell-1)(-1)^{\ell+1}}.
\end{array}
$$ 
On the other hand
$$
\begin{array}{lll}
E_n^2 &=& \displaystyle\prod_{k=1}^{2n} e^{1/2} \left(1 - \frac{1}{k+1}\right)^{(k^2+k)(-1)^k} 
    = \displaystyle\prod_{\ell=2}^{2n+1} e^{1/2} \left(1 - \frac{1}{\ell}\right)^{(\ell^2-\ell)(-1)^{\ell+1}} \\ 
   &=& e^n \displaystyle\prod_{\ell=2}^{2n+1}\left(1 - \frac{1}{\ell}\right)^{\ell^2(-1)^{\ell+1}}
           \times \displaystyle\prod_{\ell=2}^{2n+1}\left(1 - \frac{1}{\ell}\right)^{\ell(-1)^{\ell}}.
\end{array}
$$ 
Multiplying out the two expressions obtained for $E_n^2$ yields
$$
\begin{array}{lll}
E_n^4 &=& 2 e^{2n} \displaystyle\left(\frac{2n}{2n+1}\right)^{(2n+1)^2} \times
          \displaystyle\prod_{k=2}^{2n}\left(\frac{1 - \frac{1}{k+1}}{1 - \frac{1}{k}}\right)^{k^2(-1)^k} 
          \times \prod_{\ell=2}^{2n+1}\left(1 - \frac{1}{\ell}\right)^{(-1)^{\ell}} \\
      &=& 2 e^{2n} \displaystyle\left(\frac{2n}{2n+1}\right)^{(2n+1)^2} \times
          \displaystyle\prod_{k=2}^{2n}\left(\frac{k^2-1}{k^2}\right)^{-k^2(-1)^k} 
          \times \prod_{\ell=2}^{2n+1}\left(1 - \frac{1}{\ell}\right)^{(-1)^{\ell}}.
\end{array}
$$
But $2 e^{2n} \displaystyle\left(\frac{2n}{2n+1}\right)^{(2n+1)^2}$ tends to $2e^{-3/2}$ when $n$
tends to infinity (take the logarithm). We also have (see \cite[p.\,287]{Adam})
$$
\lim_{n \to \infty} \prod_{k=2}^{2n}\left(\frac{k^2-1}{k^2}\right)^{-k^2(-1)^k}
= \lim_{n \to \infty} \prod_{k=2}^{2n}\left(1-\frac{1}{k^2}\right)^{-k^2(-1)^k}
= \frac{\pi}{4}\exp\left(\frac{1}{2} + \frac{7\zeta(3)}{\pi^2}\right)
$$
and
$$
\begin{array}{lll}
\displaystyle\lim_{n \to \infty}\prod_{\ell=2}^{2n+1}\left(1 - \frac{1}{\ell}\right)^{(-1)^{\ell}}
&=& \displaystyle\lim_{n \to \infty}\left(\prod_{k=1}^n\left(1 - \frac{1}{2k}\right)
    \displaystyle\prod_{k=1}^n\left(1 - \frac{1}{2k+1}\right)^{-1}\right) \\
&=& \displaystyle\lim_{n \to \infty}\prod_{k=1}^n\left(\frac{(2k-1)(2k+1)}{(2k)^2}\right)
= \displaystyle\lim_{n \to \infty}\prod_{k=0}^{n-1}\left(\frac{(2k+1)(2k+3)}{(2k+2)^2}\right) \\
&=& \displaystyle\lim_{n \to \infty}\prod_{k=0}^{n-1}\left(\frac{(k+1/2)(k+3/2)}{(k+1)^2}\right) 
= \displaystyle\frac{\Gamma(1)^2}{\Gamma\left(\frac{1}{2}\right)\Gamma\left(\frac{3}{2}\right)}
= \displaystyle\frac{2}{\Gamma\left(\frac{1}{2}\right)^2} = \frac{2}{\pi}
\end{array}
$$
(see, e.g., \cite[Section~12-13]{WW}). Hence, finally,
$$
\lim_{n \to \infty} E_n^4 = \exp\left(-1+\frac{7\zeta(3)}{\pi^2}\right), \ \mbox{\rm thus} \
\lim_{n \to \infty} E_n = \exp\left(-\frac{1}{4} + \frac{7\zeta(3)}{4\pi^2}\right)
$$
which is Formula~(\ref{2}). \endpf

\begin{remark}
In \cite{KT} the authors note that multiplying Formulas~(\ref{1}) and (\ref{2}) together
and squaring imply the following relation
$$
\lim_{k \to \infty} \left(\frac{2^{2^2} \cdot 4^{4^2} \cdot 6^{6^2} \cdots (2k)^{(2k)^2}}
{1^{1^2} \cdot 3^{3^2} \cdot 5^{5^2} \cdots (2k-1)^{(2k-1)^2}}\right)^4
\left(\frac{(2k+2)^{4k+5}}{(2k+1)^{12k+9}}\right)^k = \exp\left(\frac{7\zeta(3)}{\pi^2}\right)
$$
which they show equivalent to the formula given by Guillera and Sondow in \cite[Example~5.3]{GS}
$$
\left(\frac{2^1}{1^1}\right)^{\frac{1 \cdot 2}{2^4}}
\left(\frac{2^2}{1^1 \cdot 3^1}\right)^{\frac{2 \cdot 3}{2^5}}
\left(\frac{2^3 \cdot 4^1}{1^1 \cdot 3^3}\right)^{\frac{3 \cdot 4}{2^6}}
\left(\frac{2^4 \cdot 4^4}{1^1 \cdot 3^6 \cdot 5^1}\right)^{\frac{4 \cdot 5}{2^7}} \cdots 
= \exp\left(\frac{7\zeta(3)}{4\pi^2}\right).
$$
The authors of \cite{KT} also note that Formula~(\ref{3}) is a rearrangement
of a formula given by Guillera and Sondow in \cite[Example~5.5]{GS}
$$
\left(\frac{3^1}{1^1}\right)^{\frac{1}{2^3}}
\left(\frac{3^2}{1^1 \cdot 5^1}\right)^{\frac{2}{2^4}}
\left(\frac{3^3 \cdot 7^1}{1^1 \cdot 5^3}\right)^{\frac{3}{2^5}}
\left(\frac{3^4 \cdot 7^4}{1^1 \cdot 5^6 \cdot 9^1}\right)^{\frac{4}{2^6}} \cdots
= \exp\left(\frac{G}{\pi}\right). 
$$
which is in turn equivalent to
$$
\lim_{k \to \infty} \left(\frac{3^3 \cdot 7^7 \cdot 11^{11} \cdots (4k-1)^{4k-1}}
{1^1 \cdot 5^5 \cdot 9^9 \cdots (4k-3)^{4k-3}}\right)^2 \frac{(4k+3)^{2k+1}}{(4k+1)^{6k+1}}
= \exp\left(\frac{4G}{\pi}\right).
$$
We thus see that both formulas in \cite[Example~5.3]{GS} and \cite[Example~5.5]{GS}
can be deduced from known values of the functions $D$ and $E$ in \cite{BD} and \cite{Adam}.
\end{remark}

\section{Conclusion}

In \cite{SH} the authors note that Borwein-Dykshoorn formulas
$$
\lim_{n \to \infty} \prod_{n=1}^{2N+1} \left(1+\frac{1}{n}\right)^{n(-1)^{n+1}} =
e \lim_{n \to \infty} \prod_{n=1}^{2N} \left(1+\frac{1}{n}\right)^{n(-1)^{n+1}} =
\frac{A^6}{2^{1/6}\sqrt{\pi}}
$$
can be written
$$
\prod_{n=1}^{\infty} \left(\frac{e}{\left(1+\frac{1}{n}\right)^n}\right)^{(-1)^{n-1}} = 
\frac{2^{1/6} e \sqrt{\pi}}{A^6}.
$$
A similar reasoning proves that
$$
D(x) = \lim_{k \to \infty} \prod_{k=1}^{2n+1} \left(1+\frac{x}{k}\right)^{k(-1)^{k+1}}
= e^x \prod_{k=1}^{\infty}\left(e^{-x}\left(1+\frac{x}{k}\right)^k\right)^{(-1)^{k+1}}.
$$
This in turn implies that
$$
\log D(x) = x + \sum_{k=1}^{\infty} (-1)^{k+1} \left(-x+k\log\left(1+\frac{x}{k}\right)\right).
$$ 
Now recall the definition of the ``parameterized-Euler-constant function'' $\gamma_{\alpha}(z)$
defined in \cite[Definition~3.1]{Xia} for $|z| \leq 1$ and $\alpha > -1$ by
$$
\gamma_{\alpha}(z) = 
\sum_{n=1}^{\infty} z^{n-1} \left(\frac{\alpha}{n} - \log \left(1+\frac{\alpha}{n}\right)\right).
$$
For $|z| < 1$ we have
$$
\gamma_{\alpha}(z) + z \gamma'_{\alpha}(z) = 
\sum_{n=1}^{\infty} z^{n-1} \left(\alpha - n \log \left(1+\frac{\alpha}{n}\right)\right).
$$ 
Thus (with the same justification as in the proof of \cite[Theorem~16]{SH}) we have the
following relation between $D$ and $\gamma_x$
$$
D(x) = e^{1+ \gamma'_x(-1) - \gamma_x(-1)}.
$$

After having read a first version of this paper on ArXiv, K. Hessami Pilehrood indicated to us
that Xia's function is actually a particular case of the function $\gamma_{a,b}(z)$ 
introduced and studied in \cite{HH}
$$
\gamma_{a,b}(z) = 
\sum_{n=0}^{\infty} \left(\frac{1}{an+b} - \log\left(\frac{an+b+1}{an+b}\right)\right)z^n.
$$
This definition is \cite[Relation~(14)]{HH} ($a, b$ positive integers, $|z| \leq 1$), while
\cite[Theorem~1]{HH} gives the analytic continuation of $\gamma_{a,b}(z)$ for
$a, b$ positive reals and $z \in {\mathbb C} \setminus [1, +\infty)$. It is clear that
$$
\gamma_{1/\alpha,1/\alpha}(z) = \gamma_{\alpha}(z)
$$
(the function on the left side is the one in \cite{HH}, the one on the right side is the one 
in \cite{Xia}).

\smallskip

\noindent
Note that, in view of \cite[Corollary~3]{HH} (see also \cite[3.6]{Xia}), this gives an 
expression of $D(x)$ in terms of the Lerch transcendent (see \cite{GS, SH}) $\Phi(z,s,u) =
\sum_{n \geq 0} \frac{z^n}{(n+u)^s}$ and its derivatives. It is then no real surprise that 
the quantities $7\zeta(3)/4\pi^2$ and $G/\pi$ also occur in Examples~2.2 and 2.3 of \cite{GS}
in the relations
$$
\frac{\partial \Phi}{\partial s}(-1,-2,1) = \frac{7\zeta(3)}{4\pi^2} \ \ \mbox{\rm and} \ \
\frac{\partial \Phi}{\partial s}(-1,-1,\frac{1}{2}) = \frac{G}{\pi}\cdot
$$

\begin{remark} It is worth noting that Equation~\ref{3} (hence also Equation~\ref{4}) can be proved 
directly from the paper of Adamchik \cite{Adam} by using a result of Choi and Srivastava \cite{CS}. 
Namely recall that the Barnes function $G(z)$ is defined by $G(1) = 1$ and $G(z+1) = G(z) \Gamma(z)$. 
Adamchik proved in \cite[Proposition~5, p.\ 284]{Adam} the following equality for $\Re(x) > -1/2$:
$$
\lim_{n \to \infty} \prod_{k=1}^{2n} \left(1 + \frac{2x}{k}\right)^{-k(-1)^k} =
\frac{e^{-x}\,\Gamma(x+\frac{1}{2})}{\Gamma(\frac{1}{2})}
\left(\frac{G(x+\frac{1}{2})}{G(x+1)G(\frac{1}{2})}\right)^2.
$$
Putting succesively $x = + 1/4$ and $x = -1/4$ and taking the quotient of the two resulting limits yields
$$
\begin{array}{lll}
\displaystyle\lim_{n \to \infty} \prod_{k=1}^{2n} \frac{\left(1 + \frac{1}{2k}\right)^{-k(-1)^k}}
{\left(1 - \frac{1}{2k}\right)^{-k(-1)^k}} 
&=& e^{-1/2} \displaystyle\frac{\Gamma(\frac{3}{4})}{\Gamma(\frac{1}{4})}
\left(\frac{G(\frac{3}{4})^2}{G(\frac{5}{4})G(\frac{1}{4})}\right)^2 \\
&=& e^{-1/2} \displaystyle\frac{\Gamma(\frac{3}{4})}{(\Gamma(\frac{1}{4}))^3}
\left(\frac{G(\frac{3}{4})}{G(\frac{1}{4})}\right)^4 \ 
\mbox{\rm since $G(\frac{5}{4}) = G(1 + \frac{1}{4}) = G(\frac{1}{4}) \Gamma(\frac{1}{4})$} \\
&=& e^{-1/2} \displaystyle\Gamma(\frac{3}{4})\Gamma(\frac{1}{4})
\left(\frac{G(\frac{3}{4})}{G(\frac{1}{4})\Gamma(\frac{1}{4})}\right)^4
\end{array}
$$
But we have from \cite[(1.15) p.\ 94]{CS}
$$
\frac{G(\frac{3}{4})}{G(\frac{1}{4})\Gamma(\frac{1}{4})} = 2^{-1/8} \pi^{-1/4} e^{G/2\pi}.
$$
Hence
$$
\lim_{n \to \infty} \prod_{k=1}^{2n} \frac{\left(1 + \frac{1}{2k}\right)^{-k(-1)^k}} 
{\left(1 - \frac{1}{2k}\right)^{-k(-1)^k}} =
\frac{e^{-1/2} \Gamma(\frac{3}{4}) \Gamma(\frac{1}{4}) e^{2G/\pi}}{\pi \sqrt{2}}
= e^{-1/2} e^{2G/\pi} 
$$
(using as above that $\Gamma(\frac{3}{4}) \Gamma(\frac{1}{4}) = \pi \sqrt{2}$). 
This is clearly equivalent to Equation~\ref{3}.

\end{remark}

\begin{remark}
A product resembling the products studied in \cite{KT} is given by Holcombe in \cite{Holc}:
$$
\pi = e^{3/2} \prod_{n=2}^{\infty} e \left(1 - \frac{1}{n^2}\right)^{n^2}.
$$
\end{remark}

\medskip

\noindent
{\bf Acknowledgements.} We would like to thank Khodabakhsh Hessami Pilehrood, Jia-Yan Yao, and the referee 
for their useful comments on a previous version of this paper.

\end{document}